\documentclass[12pt]{article}
\begin{document}
\marginpar{\tiny March 22/05}

\begin{center}
{\Large Two models of partial differential equations with discrete
and distributed state-dependent delays}

\bigskip

{\sc Alexander V. Rezounenko}

\smallskip

Department of Mechanics and Mathematics, Kharkov University,

 4, Svobody Sqr., Kharkov, 61077, Ukraine

 E-mail : rezounenko@univer.kharkov.ua

\end{center}

\begin{quote}
{\bf Abstract.} This work is the first attempt to treat partial
differential equations with discrete (concentrated)
state-dependent delay. The main idea is to approximate the
discrete delay term  by a sequence of distributed delay terms (all
with state-dependent delays). We study local existence and
long-time asymptotic behavior of solutions and prove that the
model with distributed delay has a global attractor while the one
with discrete delay possesses the trajectory attractor.
\end{quote}

\medskip



\bigskip
{\bf 1. Introduction} 

The theory of delay ordinary differential equations has a rich
history and still be one of the actively developing branches of
the theory of differential equations. We cite a few monographs
which are the classical source of fundamental facts and approaches
in this field \cite{Hale,Mishkis,Hale_book,Walther_book,Azbelev}.

Another developed branch of the theory of differential equations
is the theory of partial differential equations (PDEs). We refer
to \cite{Babin-Vishik,Temam_book,Chueshov_book}  where many deep
results on the qualitative theory of PDEs are presented.

These fields have very much in common when we are interested in
the qualitative theory. It is not surprisingly since both delay
equations and PDEs can be treated as abstract dynamical systems in
infinite-dimensional spaces. Recently, some efforts have been
applied to develop the theory of PDEs with delay. Such equations
are naturally more difficult since they are infinite-dimensional
in both time and space variables. We refer to the monograph
\cite{Wu_book} and to a few articles which are close to the
subject of this work \cite{travis_webb,Chueshov-JSM-1992,
Cras-1995,NA-1998,Rezounenko-Wu-2005}.

Recently, the theory of state-dependent ordinary differential
equations (equations where delay depends on the state of the
system) has attracted attention of many researches. We refer to
\cite{Nussbaum-Mallet-1992,Nussbaum-Mallet-1996,MalletParet,Walther2002,
Walther_JDE-2003} and references therein. The approach in these
works essentially based on the Lipschitz continuity in time of
solutions of ordinary differential equations. Unfortunately, the
last property does {\it not} hold for solutions of PDEs, so one
has to propose a new approach.

The first attempt to treat PDEs with state-dependent delay has
been done in \cite{Rezounenko-Wu-2005}. There was proposed a model
of PDEs with distributed state-dependent (state-selective) delay;
the existence and uniqueness of solutions have been proved and the
asymptotic behavior of solutions has been studied.

The present article is the first attempt to treat PDEs with
discrete (state-dependent) delay. We propose two models of PDEs
with discrete and distributed (state-dependent) delays and study
their local and long-time asymptotic behavior.
 The main idea of the present work
is to approximate the discrete delay term  by a sequence of
distributed delay terms (cf. the forms of $F$ and $F_n$ in
(\ref{sdd3-1}) and (\ref{sdd3-g})). We first develop the
techniques for studying PDEs with distributed (state-dependent)
delay and than apply it to investigations of PDEs with discrete
(state-dependent) delay. We propose a sequence of simple
distributed delay terms constructed as integrals over $(-r,0)$
with step functions as kernels of these integrals. More precisely,
using the well-known Lebesgue theorem we approximate the value of
$y(s)$ for almost all $s\in (a,b)$ by the sequence $\left\{
\varepsilon^{-1}_n\int^s_{s-\varepsilon_n} y(\tau)d\tau
\right\}^\infty_{n=1},$ where $\varepsilon_n\to 0_{+}$ as
$n\to\infty.$ These integrals can be rewritten in the form
$\int^0_{-r} y(s+\theta)\cdot\widetilde{\xi}^n(\theta,s)d\theta,$
where $\widetilde{\xi}^n$ is the step-function
$\widetilde{\xi}^n(\theta,s)\equiv\varepsilon^{-1}_n$ for
$\theta\in [s-\varepsilon_n,s]$ and
$\widetilde{\xi}^n(\theta,s)\equiv 0$ for $\theta\notin
[s-\varepsilon_n,s].$

For the model with distributed delay we prove (section~3) the
existence and uniqueness theorems, construct an evolution
semigroup and obtain the existence of global attractor. Since for
the model with discrete state-dependent delay the uniqueness of
solutions is not assumed, to study the long-time asymptotic
dynamics of these solutions we apply (section~4) the theory of
trajectory attractors (see \cite{Chepyzhov-Vishik_JMPA-1997} and
references therein).

The obtained results can be applied to the diffusive Nicholson's
blowflies equation (see e.g. So and Yang (1998), So, Wu and Yang
(2000)) with state-dependent (both discrete and distributed)
delays.

\medskip
{\bf 2. Formulation of the models with discrete and distributed
delays}

\medskip

Let us start with the following non-local partial differential
equation with {\it state-dependent \underline{discrete} delay}
\begin{equation}\label{sdd3-1}
\begin{array}{lll}
&\,\,\,\, \frac{\partial }{\partial t}u(t,x)+Au(t,x)+du(t,x)\\
&= 
\int_\Omega b(u(t-\eta (u(t),u_t), y)) f(x-y) dy   \equiv \big(
F(u_t) \big)(x),\quad x\in \Omega ,\end{array}
\end{equation}
 where $A$ is a densely-defined self-adjoint positive linear operator
 with domain $D(A)\subset L^2(\Omega )$ and with compact
  resolvent, so $A: D(A)\to L^2(\Omega )$ generates an analytic semigroup,
  $\Omega $ is a smooth bounded domain in $R^{n_0}$, $f: \Omega -\Omega \to R$
  is a bounded function to
  be specified later, $b:R\to R$ is a locally Lipschitz bounded map
  ($|b(w)|\le C_b$ with $C_b\ge 0),$ $d$ is a positive constant. The function
  $\eta (\cdot,\cdot): L^2(\Omega)\times L^2(-r,0;L^2(\Omega)) \to R$ represents the
state-dependent {\it discrete} delay. We denote for short $H\equiv
L^2(\Omega)\times L^2(-r,0;L^2(\Omega)).$

Consider the following non-local partial differential equation
with {\it state-dependent \underline{distributed} delay}
\begin{equation}\label{sdd3-g}
\begin{array}{lll}
&\quad \frac{\partial }{\partial t}u(t,x)+Au(t,x)+du(t,x)\\
&= 
\int^0_{-r} \left\{ \int_\Omega b(u(t+\theta, y))
f(x-y) dy \right\} \xi^n(\theta , u(t), u_t) d\theta \\
&\equiv \big( F_n(u_t) \big)(x), x\in \Omega ,\end{array}
\end{equation}
 where the function $\xi^n (\cdot,\cdot,\cdot): [-r,0]\times
H
 \to R$ represents the
state-dependent {\it distributed} delay.

We consider equations (\ref{sdd3-1}) or (\ref{sdd3-g}) with the
following initial conditions
\begin{equation}\label{sdd3-ic}
  u(0+)=u^0\in L^2(\Omega), \quad  u|_{(-r,0)}=\varphi \in L^2 (0,T;L^2(\Omega)).
\end{equation}

\medskip
{\bf 3. Distributed delay problem} \medskip

In this section we study the existence and properties of solutions
for distributed delay problem (\ref{sdd3-g}), (\ref{sdd3-ic}).

\medskip

\noindent {\bf Definition 1.} {\it A function $u$ is a {\it weak
solution} of problem (\ref{sdd3-g}) subject to the initial
conditions (\ref{sdd3-ic}) on an interval $[0,T]$ if $u\in
L^\infty (0,T;L^2(\Omega))\cap
 L^2 (-r,T;L^2(\Omega))\cap  L^2 (0,T;D(A^{1\over 2}))$, $u(\theta)=\varphi (\theta)$ for $\theta\in (-r,0)$
 and
\begin{equation}\label{sdd3-sol}
  -\int^T_0\langle u, \dot v\rangle dt +
 \int^T_0\langle A^{1\over 2}u,A^{1\over 2}v\rangle dt +
\int^T_0\langle du-F_n(u_t),v\rangle dt = -\langle u^0,v(0)\rangle
\end{equation}
 for any function $v\in L^2 (0,T;D(A^{1\over 2}))$ with
 $\dot v\in L^2 (0,T;D(A^{-{1\over 2}}))$ and $v(T)=0.$

 }

\medskip
\noindent {\bf Theorem 1.} {\it Assume that
\begin{itemize}
\item[(i)] $b: R\to R$ is locally Lipschitz and bounded i.e.,
there
exists a constant $C_b$ so that that $|b(w)|\le C_b$ 
for all $w\in R$; \item[(ii)] $f: \Omega-\Omega\to R$ is bounded;
\item[(iii)] $\xi^n : [-r,0]\times L^2(\Omega)\times
L^2(-r,0;L^2(\Omega))\to R$ satisfies the
following conditions:\\
{\bf a)} for any $M>0$ there exists $L_{\xi , M,n}$ so that for
all $(v^i, \psi^i)\in H$ satisfying
 $||v^i||^2+ \int^0_{-r} ||\psi^i(s)||^2 ds\le M^2, i=1,2$
 one has
\begin{equation}\label{sdd3-xi}
\quad \int^0_{-r}|\xi^n (\theta, v^1, \psi^1)-\xi^n (\theta, v^2, \psi^2) | d\theta\\
\le L_{\xi,M,n}
  \cdot  \left( ||v^1-v^2||^2+ \int^0_{-r} ||\psi^1(s)-\psi^2(s)||^2ds
  \right)^{1/2},
\end{equation}
{\bf b)} there exists $C_{\xi,1} >0$ so that
\begin{equation}\label{sdd3-cxi}
\Vert \xi^n(\cdot\, , v, \psi)\Vert_{L^1(-r,0)} \le C_{\xi,1}
\,\,\mbox{for all}\,\,  (v, \psi)\in H.
\end{equation}
\end{itemize}

Then for any $(u^0, \varphi) \in H\equiv L^2(\Omega)\times L^2
(-r,0;L^2(\Omega))$ the problem (\ref{sdd3-g}) subject to the
initial conditions (\ref{sdd3-ic}) has a weak solution $u(t)$ on
every given interval $[0,T]$
 which satisfies
\begin{equation}\label{sdd3-contin} u(t)\in C([0,T];L^2(\Omega)).
\end{equation}
}

\medskip
\noindent {\it Proof of Theorem~1}. Let us denote by $\{
e_k\}^\infty_{k=1}$ an orthonormal basis of $L^2(\Omega)$ such
that $Ae_k=\lambda_ke_k$, $0< \lambda_1<\ldots<\lambda_k\to
+\infty$. We say that function
$u^m(t,x)=\sum\limits^m_{k=1}g_{k,m}(t)e_k(x) $ is a {\it Galerkin
approximate solution of order $m$ for the problem
(\ref{sdd3-g}),(\ref{sdd3-ic})} if
\begin{equation}\label{sdd3-6}
\left\{ \begin{array}{ll} &\langle \dot u^m+Au^m
+du^m-F_n(u^m_{t}),
e_k\rangle =0,\\
&\langle u^m(0+),e_k\rangle =\langle u^0,e_k\rangle , \,\, \langle
u^m(\theta),e_k\rangle=\langle \varphi (\theta) , e_k\rangle
,\,\,\forall \theta\in (-r,0) \end{array}\right.
 \end{equation}
$\forall k=1,\ldots,m$. Here $g_{k,m}\in C^1(0,T;R)\cap
L^2(-r,T;R)$ with $\dot g_{k,m}(t)$ being absolutely continuous.

Equations (\ref{sdd3-6}) for fixed $m$ and $n$ can be rewritten as
the following system for the $m$-dimensional vector-function
$v(t)=v^m(t)=(g_{1,m}(t),\ldots,g_{m,m}(t))^T :$
\begin{equation}\label{sdd3-7}
\dot v(t)=\hat f(v(t))+ \int^0_{-r} p(v(t+\theta))
\widetilde\xi^n(\theta,v(t),v_t) d\theta,
\end{equation}
where function $ \widetilde\xi^n$ satisfies properties similar to
(\ref{sdd3-xi}), (\ref{sdd3-cxi}) if one uses $|\cdot|_{R^m}$
instead of $\Vert \cdot\Vert_{L^2(\Omega)}.$ We notice that $\Vert
u^m(t,\cdot)\Vert^2_{L^2(\Omega)}=
 \sum\limits^m_{k=1}g^2_{k,m}(t)=|v(t)|^2_{R^m}.$

Under the above assumptions, the functions $\hat f$ and $p$ are
locally Lipschitz, $|p(s)|\le c_2$ for $s\in R$,
Therefore, for any initial data $\varphi\in L^2(-r,0; R^m),$ $a\in
R^m$ Theorem~6 and Remark~9 from Rezounenko (2004) give that there
exists $\alpha>0$ and a unique solution of (\ref{sdd3-7}) $v\in
L^2(-r,\alpha; R^m)$ such that $v_0=\varphi $ and $v(0)=a$, and
$v|_{[0,\alpha]}\in C([0,\alpha]; R^m)$ (for more details see
Rezounenko (2004)).

It is easy to get from the boundedness of $b$ and (\ref{sdd3-cxi})
that
\begin{equation}\label{sdd3-F1}
 |\langle F_n(u_t),v\rangle_{L^2(\Omega)}|
\le  M_f|\Omega |^{3/2} C_b C_{\xi,1}\cdot \Vert v\Vert.
\end{equation}

Now, we try to get an {\it a-priori} estimate for the Galerkin
approximate solutions for the problem
(\ref{sdd3-g}),(\ref{sdd3-ic}). We multiply (\ref{sdd3-6}) by
$g_{k,m}$ and sum over $k=1,\cdots ,m$. Hence for  $u(t)=u^m(t)$
and $t\in (0,\alpha]\equiv(0,\alpha (m)]$, the local existence
interval for $u^m(t)$, we get
\begin{equation}\label{sdd3-11}
{1\over 2}{d\over dt}\Vert u(t)\Vert^2 + \Vert A^{1/2}u(t)\Vert^2
+ d \Vert u(t)\Vert^2  \le |\langle F_n(u_t),u(t)\rangle|.
\end{equation}
Using (\ref{sdd3-F1}), we obtain
\begin{equation}\label{sdd3-12}
{d\over dt}\Vert u(t)\Vert^2 + 2\Vert A^{1/2}u(t)\Vert^2 \le
\tilde k_1 \Vert u(t)\Vert^2
+\tilde k_3.
\end{equation}
Since ${d\over dt}\Vert u(t)\Vert^2 + 2\Vert A^{1/2}u(t)\Vert^2 =
{d\over dt}\left(\Vert u(t)\Vert^2 + 2\int^t_0\Vert
A^{1/2}u(\tau)\Vert^2d\tau +\tilde k_3\right) $, we denote by
$\chi (t)\equiv \Vert u(t)\Vert^2 + 2\int^t_0\Vert
A^{1/2}u(\tau)\Vert^2d\tau +\tilde k_3$ and rewrite the last
estimate as follows $ {d\over dt}\chi (t) \le
\tilde k_1\cdot\chi (t).$ Multiplying it by $e^{-\tilde k_1 t},$
one gets $ {d\over dt}\left(e^{-\tilde k_1t}\chi (t) \right) \le
0.$
Integrating from $0$ to $t$ and then
multiplying by $e^{\tilde k_1 t},$  we obtain 
$ \chi(t)\le \left( \Vert u(0)\Vert^2 +\tilde k_3\right) e^{\tilde
k_1 t}. $ So, we have the {\it a -priori} estimate
\begin{equation}\label{sdd3-14}
\Vert u(t)\Vert^2 + 2\int^t_0\Vert A^{1/2}u(\tau)\Vert^2d\tau\le
\left( \Vert u(0)\Vert^2+\tilde k_3\right)e^{\tilde k_1 t} -
\tilde k_3.
\end{equation}
 Estimate (\ref{sdd3-14}) gives that, for
$u^0\in L^2(\Omega)$ 
the family of approximate solutions $\{ u^m(t)\}^\infty_{m=1}$ is
uniformly (with respect to $m\in {\bf N}$) bounded in the space
$L^\infty (0,T;L^2(\Omega))\cap L^2 (0,T;D(A^{1/2})),$ where
$D(A^{1/2})$ is the domain of the operator $A^{1/2}$ and $[0,T]$
is the local existence interval. From (\ref{sdd3-14}) we also get
the continuation of $u^m(t)$ on any interval, so (\ref{sdd3-14})
holds for all $t>0.$

Using the definition of Galerkin approximate solutions
(\ref{sdd3-6}) and their property (\ref{sdd3-14}), we can
integrate over $[0,T]$ to obtain $\int^T_0\Vert A^{-{1/2}} \dot
u^m(\tau)\Vert^2 d\tau \le C_T$ for any $T.$ These properties of
the family $\{ u^m(t)\}^\infty_{m=1}$ give that $\{ (u^m(t); \dot
u^m(t))
\}^\infty_{m=1}$ is a bounded sequence in the space 
\begin{equation}\label{sdd3-13}
X_T\equiv L^\infty (0,T;L^2(\Omega))\cap L^2
(0,T;D(A^{1/2}))\times L^2(0,T;D(A^{-{1/2}})).
\end{equation}
 Then
there exist a function $(u(t); \dot u(t))$ and a subsequence $\{
u^{m_k} \} \subset \{ u^m \}$ such that
\begin{equation}\label{sdd3-*}
(u^{m_k}; \dot u^{m_k}) \quad \hbox{*-weakly converges to}
 \quad  (u; \dot u) \quad \hbox{in the space} \quad X_T.
\end{equation}
By a standard argument (using the strong convergence $u^{m_k}\to
u$ in the space $L^2 (0,T;L^2(\Omega))$ which follows from
(\ref{sdd3-*}) and the Doubinskii's theorem, one can show (see
e.g. Lions (1969), Chueshov (1999) and Rezounenko (1997)) that any
*-weak limit is a solution of (\ref{sdd3-g}) subject to the
initial conditions (\ref{sdd3-ic}). To prove the continuity of
weak solutions we use the well-known 

\medskip
\noindent {\bf Proposition~1} (see e.g., Proposition 1.2 in
\cite{showalter}). {\it Let the Banach space $V$ be dense and
continuously embedded in the Hilbert space $X;$ identify $X=X^*$
so that $V\hookrightarrow X\hookrightarrow V^*.$ Then the Banach
space $W_p(0,T)\equiv \{ u\in L^p(0,T;V) : \dot u\in
L^q(0,T;V^*)\}$ (here $p^{-1}+q^{-1}=1$) is contained in
$C([0,T];X).$ }

In our case $X=L^2(\Omega), V=D(A^{{1/2}}), V^*=D(A^{-{1/2}}),
p=q=1/2$ (see (\ref{sdd3-13}),(\ref{sdd3-*})). Hence Proposition~1
gives (\ref{sdd3-contin}). The proof of Theorem~1 is complete.
\rule{5pt}{5pt}

\bigskip

Now we describe a sufficient condition for the uniqueness of weak
solutions.

\noindent {\bf Theorem~2.} {\it Assume that functions $b$ and $f$
are as in Theorem~1 (satisfy properties (i),(ii)), function
$\xi^n$ satisfies property (iii)-a) and
\begin{equation}\label{sdd3-23}
   \xi^n (\cdot, v,\psi) \in L^\infty(-r,0)\,\,\mbox{for all}\,\,  (v, \psi)\in H.
\end{equation}}
Then solution of (\ref{sdd3-g}), (\ref{sdd3-ic}) given by
Theorem~1  is unique.

\medskip

\noindent {\it Proof of Theorem~2.} Let $u^1$ and $u^2$ be two
solutions of (\ref{sdd3-g}), (\ref{sdd3-ic}).
Below we denote for short
$w(t)=w^{n,m}(t)=u^{1,n,m}(t)-u^{2,n,m}(t)$ - the difference of
corresponding Galerkin approximate solutions. Hence
\begin{equation}\label{sdd3-24}
{d\over dt}\Vert w(t)\Vert^2 + 2 \Vert A^{1/2}w(t)\Vert^2 + 2d
\Vert w(t)\Vert^2 = \langle F_n(u^1_t)- F_n(u^{2}_t),w(t)\rangle.
\end{equation}
Let us consider the difference $\langle F_n(u^1_t)-
F_n(u^{2}_t),w(t)\rangle$ in details  (see (\ref{sdd3-1}),
(\ref{sdd3-g})).
$$\langle F_n(u^1_t)-
F_n(u^{2}_t),w(t)\rangle\equiv \int_\Omega\left[
 \int^0_{-r} \left\{ \int_\Omega b(u^1(t+\theta, y)) f(x-y) dy
\right\} \xi^n(\theta , u^1(t), u^1_t) d\theta  - \right.$$
$$
-\left. \int^0_{-r} \left\{ \int_\Omega b(u^{2}(t+\theta, y))
f(x-y) dy \right\} \xi^n(\theta , u^{2}(t), u^{2}_t) d\theta
\right] \cdot w(t,x) dx
$$
$$=\int_\Omega\left[ \int^0_{-r} \left\{ \int_\Omega b(u^1(t+\theta,
y)) f(x-y) dy \right\} \xi^n(\theta , u^1(t), u^1_t) d\theta
-\right.
$$
$$
-\left. \int^0_{-r} \left\{ \int_\Omega b(u^{2}(t+\theta, y))
f(x-y) dy \right\} \xi^n(\theta , u^1(t), u^1_t) d\theta \right]
\cdot w(t,x) dx,
$$
$$
+\int_\Omega\left[ \int^0_{-r} \left\{ \int_\Omega
b(u^{2}(t+\theta, y)) f(x-y) dy \right\} \xi^n(\theta , u^1(t),
u^1_t) d\theta -\right.
$$
$$
-\left. \int^0_{-r} \left\{ \int_\Omega b(u^{2}(t+\theta, y))
f(x-y) dy \right\} \xi^n(\theta,u^{2}(t),u^{2}_t) d\theta \right]
\cdot w(t,x) dx.
$$
Using the local Lipschitz property of $b,$ (\ref{sdd3-23}) and
(\ref{sdd3-xi}), one easily checks that there are positive
constants $C_3,C_4$ such that
$$|\langle F_n(u^1_t)-
F_n(u^{2}_t),w(t)\rangle | \le C_3||w(t)||^2 + C_4\int^0_{-r}
||w(t+\theta)||^2 d\theta $$
$$\le C_3||w(t)||^2 + C_4\int^t_{-r}
||w(s)||^2 ds = C_3||w(t)||^2 + C_4\left(
\int^0_{-r}||w(\theta)||^2 d\theta +\int^t_{0}||w(s)||^2
ds\right).
$$
The last estimate, (\ref{sdd3-24}) and $\Vert A^{1/2}v\Vert^2\ge
\lambda_1\Vert v\Vert^2$ give
$${d\over dt}\Vert w(t)\Vert^2 + 2(\lambda_1+d)\Vert w(t)\Vert^2
\le C_3||w(t)||^2 + C_4\left( \int^0_{-r}||w(\theta)||^2 d\theta
+\int^t_{0}||w(s)||^2 ds\right).
$$
We rewrite this as
$${d\over dt}\left[ \Vert w(t)\Vert^2 + 2(\lambda_1+d)\int^t_{0}\Vert
w(s)\Vert^2 ds\right] \le C_3||w(t)||^2 + C_4\left(
\int^0_{-r}||w(\theta)||^2 d\theta +\int^t_{0}||w(s)||^2
ds\right).
$$
Hence there exists $C_5>0,$ such that for $Z(t)\equiv \Vert
w(t)\Vert^2 + 2(\lambda_1+d)\int^t_{0}\Vert w(s)\Vert^2 ds,$ we
have
$${d\over dt} Z(t) \le C_5Z(t)+ C_4\int^0_{-r}||w(\theta)||^2 d\theta.$$
Gronwall lemma implies
\begin{equation}\label{sdd3-25}
Z(t)\le \left( \Vert w(0)\Vert^2 +
C_4C^{-1}_5\int^0_{-r}||w(\theta)||^2 d\theta\right) \cdot
e^{C_5t}.
\end{equation}
The last estimate allows one to apply the well-known

\smallskip

\noindent {\bf Proposition~2.} \cite[Theorem~9]{yosida}
 {\it Let $X$ be a Banach space. Then any *-weak convergent sequence
 $\{ w_k\}^\infty_{n=1}\in X^{*}$  *-weak converges to an element
 $w_\infty\in X^{*}$ and $\Vert w_\infty\Vert_X \le\liminf_{n\to\infty} \Vert w_n\Vert_X.$
}

Hence, for the difference $u^1(t)-u^2(t)$ of two solutions we have
$$\Vert u^1(t)-u^2(t)\Vert^2 + 2(\lambda_1+d)\int^t_{0}\Vert
u^1(s)-u^2(s)\Vert^2ds $$
\begin{equation}\label{sdd3-26}
\le \left( \Vert u^1(0)-u^2(0)\Vert^2 +
C_4C^{-1}_5\int^0_{-r}||\varphi^1(\theta)-\varphi^2(\theta)||^2
d\theta\right) \cdot e^{C_5t}.
\end{equation}

We notice that by (\ref{sdd3-contin}) the difference $\Vert
u^1(t)-u^2(t)\Vert$ makes sense for all $t\in [0,T],\, \forall
T>0.$ The last estimate gives the uniqueness of solutions and
completes the proof of Theorem~2. \rule{5pt}{5pt}

\medskip

Theorems~1 and 2 allow us to define the evolution semigroup $S_t :
H\to H$, with  $H\equiv L^2(\Omega)\times L^2 (-r,0;L^2(\Omega))$,
by the formula $S_t(u^0;\varphi)\equiv (u(t);u(t+\theta)), \,
\theta\in (-r,0),$ where $u(t)$ is the weak solution of
(\ref{sdd3-g}),(\ref{sdd3-ic}). The continuity of the semigroup
with respect to time follows from (\ref{sdd3-contin}), and with
respect to initial conditions from (\ref{sdd3-26}).

For the study of long-time asymptotic properties of the above
evolution semigroup we recall (see e.g. Babin and Vishik (1969),
Temam (1988))

\medskip

\noindent {\bf Definition~2.} {\it A global attractor of the
semigroup $S_t$ is a closed bounded set ${\mathcal U}$ in $H,$
strictly invariant ($S_t{\mathcal U}={\mathcal U}$ for any $t\ge
0$), such that for any bounded set $B\subset H$ we have
$\lim\limits_{t\to +\infty} \sup \{ dist_H (S_ty, {\mathcal U}),
y\in B\} =0.$ }

As in \cite{Babin-Vishik,
Chueshov-JSM-1992,Chueshov_book,Rezounenko-MAG-1997} (see also
\cite{Rezounenko-Wu-2005}) we prove

\medskip

\noindent {\bf Theorem 3.} {\it Under the assumptions of
Theorems~1 and 2,
 the dynamical system
$(S_t;H)$ has a compact global attractor ${\mathcal U}$ which is a
bounded set in the space $H_1\equiv D(A^{\alpha})\times W$, where
$W=\{ \varphi : \varphi \in L^\infty (-r,0;D(A^{\alpha})),
\dot\varphi \in L^\infty (-r,0;D(A^{\alpha-1}))\}$,  $\alpha\le
{1\over 2}.$}

\bigskip
{\bf 4. Discrete delay problem}\nopagebreak 

\medskip

Consider the function $\eta : H \to R$ (as before
 $H\equiv L^2(\Omega)\times L^2 (-r,0;L^2(\Omega))$) which
represents the state-dependent {\it discrete} delay in equation
(\ref{sdd3-1}). Let us fix a positive sequence $\{
\varepsilon_n\}^\infty_{n=1} \subset R_{+}$ such that
$\varepsilon_n\to 0_{+}$ and define the sequence of functions
$\xi^n : [-r,0]\times H\to R$ as follows
\begin{equation}\label{sdd3-17}
\xi^n(\theta,a,\varphi) \equiv \left[
\begin{array}{cc}
  1/\varepsilon_n, & \theta\in [-\eta (a,\varphi)-\varepsilon_n,-\eta (a,\varphi)]; \\
  0, & \theta\not\in [-\eta (a,\varphi)-\varepsilon_n,-\eta (a,\varphi) ], \\
\end{array}
\right. \quad \hbox{ with} \quad \varepsilon_n>0.
\end{equation}

For example, such functions can be constructed as composition
\begin{equation}\label{sdd3-27}
    \xi^n(\theta,a,\varphi) =\widetilde\xi^n(\theta,-\eta
(a,\varphi)),
\end{equation}
where $\widetilde\xi^n(\theta,s) : [-r,0]\times R \to R$ are the
step-functions (see figure below)
$$
\widetilde\xi^n(\theta,s) \equiv \left[
\begin{array}{cc}
  1/\varepsilon_n, & \theta\in [s-\varepsilon_n,s]; \\
  0, & \theta\not\in [s-\varepsilon_n,s ], \\
\end{array}
\right. \quad \hbox{ with} \quad \varepsilon_n>0.
$$
\medskip

\begin{center}
\begin{picture}(260,100)
 \put(0,20){\vector(1,0){200}}
\put(140,0){\vector(0,1){110}}

\put(4,10){$-r$}\put(0,18){\line(0,1){5}}
\put(141,10){$0$}\put(140,18){\line(0,1){5}}
\put(67,10){$s-\varepsilon_k$}\put(70,18){\line(0,1){5}}
\put(97,90){\line(0,1){7}}

\put(30,23){$s-\varepsilon_n$}
\put(100,23){$s$}\put(70,18){\line(0,1){5}}
\put(145,75){$1/\varepsilon_n$}\put(138,80){\line(1,0){5}}
\put(145,95){$1/\varepsilon_k$}\put(138,97){\line(1,0){5}}

{\thicklines \put(0,20){\line(2,0){57}}
\put(140,20){\line(-5,0){43}}%
\put(57,80){\line(1,0){40}} \put(70,98){\line(5,0){27}}
}

\multiput(70,18)(0,18){5}{\line(0,1){7}}
\multiput(57,18)(0,8){8}{\line(0,1){5}}
\multiput(97,18)(0,8){8}{\line(0,1){5}}

\put(40,0){$\theta$} \put(100,60){$\widetilde\xi^n(\theta,s)$}
\put(100,90){$\widetilde\xi^k(\theta,s)$}
\end{picture}
\end{center}
\nopagebreak
\medskip
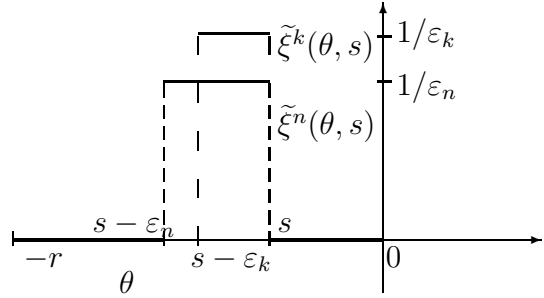
\centerline{Figure 1: Graph of functions
$\widetilde\xi^n(\theta,s)$ and $\widetilde\xi^k(\theta,s)$ with
$0<\varepsilon_k<\varepsilon_n.$ 
}

\medskip

\noindent {\bf Remark~1.} {\it It is easy to see that functions
$\xi^n,$ defined in (\ref{sdd3-17}), satisfy all the assumptions
of Theorem~2, hence weak solution of (\ref{sdd3-g}),
(\ref{sdd3-ic}) is unique for all $n.$}
\medskip

\noindent {\bf Remark~2.} {\it We notice that functions $\xi^n,$
defined in (\ref{sdd3-17}) do {\tt not} satisfy neither assumption
(iii)-a) nor assumption (iii)-b) of Theorem~1 from
\cite{Rezounenko-Wu-2005}, but do satisfy the ones of Theorem~1
given in Section~3.}

\medskip

Our idea is to use the sequence of equations (\ref{sdd3-g}) with
the right-hand sides $F_n$ (with functions $\xi^n$ defined in
(\ref{sdd3-17})) to treat equation (\ref{sdd3-1}). This idea is
based on the following well-known fact.

Let $X$ be a Banach space,  $y\in L^1(0,T;X)$ and define the
 primitive $Y(t)\equiv \int^t_0y(s)ds,\, 0\le t\le T < +\infty.$

\noindent
{\bf Proposition~3. (Lebesgue theorem)} \cite{yosida, showalter}
 {\it At a.e., $t\in (0,T), Y$ is (strongly) differentiable with
 $Y^\prime (t)\equiv \lim_{\varepsilon\to 0} {\varepsilon^{-1}}(Y(t+\varepsilon)-Y(t))=y(t).$
 }

That gives for any $(a,\varphi)\in H$
\begin{equation}\label{sdd3-28}
y(t-\eta(a,\varphi))= \lim_{n\to \infty} \int^0_{-r}
y(t+\theta)\cdot \widetilde\xi^n(\theta,-\eta(a,\varphi)) d\theta
= \lim_{n\to \infty} \int^0_{-r} y(t+\theta)\cdot
\xi^n(\theta,a,\varphi) d\theta.
\end{equation}

\medskip

\noindent {\bf Remark~3.} {\it Let us explain how to get
(\ref{sdd3-28}). If we choose in Lebesgue theorem
$\varepsilon=-\varepsilon_n<0$ and $t=s,$ one obtains
$$y(s)=\lim_{\varepsilon\to 0_{-}}
{\varepsilon^{-1}}\int^{s+\varepsilon}_s y(\theta)d\theta =
{(-\varepsilon)^{-1}}\int^s_{s+\varepsilon}
y(\theta)d\theta=\lim_{\varepsilon\to 0_{-}}
{(-\varepsilon)^{-1}}\int^0_{\varepsilon}y(s+\theta)d\theta $$
$$= \lim_{\varepsilon_n\to 0_{+}} {\varepsilon^{-1}_n}
\int^0_{-\varepsilon_n}y(s+\theta)d\theta=\lim_{\varepsilon_n\to
0_{+}} \int^0_{-r} y(s+\theta)\cdot \widetilde\xi^n(\theta,0)
d\theta$$ $$=\lim_{n\to +\infty} \int^0_{-r} y(s+h+\theta)\cdot
\widetilde\xi^n(\theta,-h) d\theta, \quad \forall h\in (0,r).$$ We
use the last equality for $h=\eta(a,\varphi)$ and
$s=t-h=t-\eta(a,\varphi)$ to get the first equality in
(\ref{sdd3-28}). The second equality in (\ref{sdd3-28}) follows
from (\ref{sdd3-27}).}

In the same way, for a function $u\in L^\infty (0,T;X)\cap
 L^2 (-r,T;X),$ one has $b(u(\cdot)) \in L^\infty (0,T;X)\cap
 L^2 (-r,T;X)\subset L^1(-r,T;X)$ and
(\ref{sdd3-28}) implies that for all $(a,\varphi)\in H$
 $$b(u(t-\eta(a,\varphi)))= \lim_{n\to \infty} \int^0_{-r}
b(u(t+\theta))\cdot \xi^n(\theta,a,\varphi) d\theta \quad \hbox{
at a.e.} \quad  t\in (-r,T).
$$
The last property reflects the main idea to approximate the
discrete delay term (with delay $\eta(a,\varphi)$) by a sequence
of distributed delay terms (cf. the forms of $F$ and $F_n$ in
(\ref{sdd3-1}) and (\ref{sdd3-g})). In Definitions~1 and 3 we are
interested in $X=L^2(\Omega).$

\medskip

\noindent {\bf Definition~3.} {\it A function $u$ is a {\tt weak
limiting solution} of problem (\ref{sdd3-1}) subject to the
initial conditions (\ref{sdd3-ic}) on an interval $[0,T]$ if $u\in
L^\infty (0,T;L^2(\Omega))\cap
 L^2 (-r,T;L^2(\Omega))\cap  L^2 (0,T;D(A^{1/ 2}))$,
 $\dot u \in L^2 (0,T;D(A^{-1/ 2}))$,
  $u(\theta)=\varphi (\theta)$ for $\theta\in (-r,0)$
 and there is a sequence $\{ (n,m) \}_{n,m\in {\bf N}},$
 such that
\begin{equation}\label{sdd3-22}
(u^{n,m}; \dot u^{n,m}) \quad \hbox{\sl
 *-weakly converges to} \quad  (u; \dot u)\quad  \hbox{in the space} \quad
 X_T, \quad \hbox{when } \min\{n,m\} \to \infty
\end{equation}
(see (\ref{sdd3-13}), (\ref{sdd3-*})). Here $u^{n,m}$ is the
Galerkin approximate solution (see (\ref{sdd3-6})) of order $m$
for the problem  (\ref{sdd3-g}), (\ref{sdd3-ic}) (the right-hand
side of (\ref{sdd3-g}) is
$F_n$ with functions $\xi^n$ defined in (\ref{sdd3-17})). 
 }

\medskip

\noindent {\bf Remark 4.} {\it It is important to note that if we
formulate in the similar manner the definition for a fixed $n$ and
condition $m\to \infty$ instead of $\min\{n,m\} \to \infty$ (see
(\ref{sdd3-22})), then we get the definition of a weak solution of
(\ref{sdd3-g}), (\ref{sdd3-ic}) which is {\tt equivalent} to
Definition~1. It follows from (\ref{sdd3-*}) and the uniqueness of
weak solutions (see Remark~1).}

\medskip

Now we prove

\smallskip

\noindent {\bf Theorem~4.} {\it Assume that functions $b$ and $f$
are as in Theorem~1 (satisfy properties (i),(ii)) and function
$\eta : H \to [0,r]$ is locally Lipschitz i.e. for any $M>0$ there
exists $L_{\eta,M}$ so that for all $(v^i, \psi^i)\in H$
satisfying
 $||v^i||^2+ \int^0_{-r} ||\psi^i(s)||^2 ds\le M^2, i=1,2$
 one has
\begin{equation}\label{sdd3-eta}
\quad |\eta (v^1, \psi^1)-\eta (v^2, \psi^2) | \\
\le L_{\eta,M}
  \cdot  \left( ||v^1-v^2||^2+ \int^0_{-r} ||\psi^1(s)-\psi^2(s)||^2ds
  \right)^{1/2},
\end{equation}

Then for any $(u^0, \varphi) \in H$ the problem (\ref{sdd3-1})
subject to the initial conditions (\ref{sdd3-ic}) has a weak
limiting solution on every given interval $[0,T]$ and
satisfies $u(t)\in C([0,T];L^2(\Omega)).$ }


\smallskip

\noindent {\it Proof of Theorem~4.} It is easy to check that
property (\ref{sdd3-eta}) implies that all functions $\xi^n,$
defined in (\ref{sdd3-17}), satisfy properties a), b) of Theorem~1
with $L_{\xi,M,n}=2\varepsilon_n^{-1}\cdot L_{\eta,M}$ and
$C_{\xi,1}=1$ for all $n.$

Now we consider any fixed sequence $\{ (n_i;m_i) \}^\infty_{i=1}$
such that $\min \{n_i;m_i\}\to \infty$ as $i\to\infty.$ Consider
the family of Galerkin approximate solutions $\{
u^{n_i,m_i}\}^\infty_{i=1}$ (see (\ref{sdd3-6})) all constructed
for the same initial data (\ref{sdd3-ic}). The a-priory estimate
(\ref{sdd3-14}) with constants $\tilde k_1, \tilde k_3$
independent of $n$ and $m$ gives that $\{ ( u^{n_i,m_i};\dot
u^{n_i,m_i}) \}^\infty_{i=1}$ is a bounded sequence in the space
$X_T$ (see (\ref{sdd3-13})). Hence there exists a *-weak
convergent subsequence, which converges (by Definition~3) to a
weak limiting solution of (\ref{sdd3-1}), (\ref{sdd3-ic}). The
continuity of a weak limiting solution follows from Proposition~1.
   The proof of Theorem~4 is complete.
\rule{5pt}{5pt}

\medskip

To study the long-time asymptotic dynamics of solutions to
(\ref{sdd3-1}), (\ref{sdd3-ic}) we apply the theory of trajectory
attractors (see \cite{Chepyzhov-Vishik_JMPA-1997} and references
therein).

\medskip

Consider the following Banach space
$$ {\cal F}^b_{+}\equiv \left\{ w(\cdot)\quad | \quad w(\cdot)\in L^b_2(R_{+}; D(A^{1/2}))\cap
L_\infty (R_{+};L^2(\Omega)),\quad  \dot  w(\cdot)\in L^b_2(R_{+};
D(A^{-{1/2}})) \right\}
$$
with the norm
$$||w||_{{\cal F}^b_{+}}= ||w||_{L^b_2(R_{+}; D(A^{1/2}))} + ||w||_{L_\infty (R_{+};L^2(\Omega))} +
||\dot w||_{L^b_2(R_{+}; D(A^{-{1/2}}))} ,
$$ where
$$||w||_{L^b_2(R_{+}; D(A^\alpha))}\equiv \sup_{h\ge 0} \int^{h+1}_h
||A^\alpha w(s)||^2 ds, \quad \hbox{ for } \alpha=1/2 \hbox{ and }
\alpha=-{1/2}.$$

Now we consider a wider space
$$ {\cal F}^{loc}_{+}\equiv \left\{ w(\cdot)\quad | \quad w(\cdot)\in
L^{loc}_2(R_{+}; D(A^{1/2}))\cap L^{loc}_\infty
(R_{+};L^2(\Omega)),\quad  \dot  w(\cdot)\in L^{loc}_2(R_{+};
D(A^{-{1/2}})) \right\}
$$
equipped with the local *-weak topology and denote it by
$\Theta^{loc}_{+}.$ More precisely, a sequence $\{ w^k \} \subset
{\cal F}^{loc}_{+}$ converges in $\Theta^{loc}_{+}$ to $ w \in
{\cal F}^{loc}_{+}$ when $k\to\infty$ if for any $T>0$
\begin{equation}\label{sdd3-21a}
(w^k; \dot w^k) \quad \hbox{\sl
 *-weakly converges to} \quad  (w; \dot w)\quad  \hbox{in the space} \quad X_T
\end{equation}
(see (\ref{sdd3-13})).  
\medskip

\noindent {\bf Remark~5.} {\it It is easy to see that ${\cal
F}^b_{+}\subset\Theta^{loc}_{+}$ and any ball $B_{R}=\{
w(\cdot)\in {\cal F}^b_{+} \quad | \quad ||w||_{{\cal F}^b_{+}}
\le R \}$ is compact in $\Theta^{loc}_{+}.$}
\medskip

\noindent {\bf Definition~4.} {\it The {\tt translation semigroup}
$\{ T(h), h\ge 0\}$ acting on the space $L^{loc}_2(R_{+};
D(A^{1/2}))\cap L^{loc}_\infty (R_{+};L^2(\Omega))$ is defined as
the set of translations along the time axis i.e.,
$$ T(h)w(\cdot)\equiv w(\cdot +h),\quad h\ge 0.
$$ }
It is evidently that the family $\{ T(h), h\ge 0\}$ is indeed a
semigroup i.e., $T(h_1+h_2)=T(h_1)T(h_2)$ for any $h_1,h_2\ge 0$
and $T(0)=Id$ -- identical operator. It is also easy to see that
the semigroup $\{ T(h), h\ge 0\}$ is {\it continuous in the
topology $\Theta^{loc}_{+}.$}

\medskip

\noindent {\bf Definition~5.} {\it {\tt Trajectory space} ${\cal
K}^{+}$ for equation (\ref{sdd3-1}) is the space of functions
$u\in {\cal F}^b_{+}$ such that for any $T>0$ the restriction
$u|_{[0,T]}$ is a weak limiting solution of (\ref{sdd3-1}),
(\ref{sdd3-ic}).}
\medskip

\noindent {\bf Remark~6.} {\it It is easy to see that trajectory
space ${\cal K}^{+}$ is invariant under the translation semigroup
$\{ T(h), h\ge 0\}$  i.e., $T(h){\cal K}^{+}\subset {\cal K}^{+},
\, \forall h\ge 0.$}

\medskip

\noindent {\bf Remark~7.} {\it By Definition~3, it is not too hard
to prove that trajectory space ${\cal K}^{+}$ is closed in the
topology $\Theta^{loc}_{+}.$}

\medskip

\noindent {\bf Definition~6.} {\it A set $P\subset {\cal
F}^{loc}_{+}$ is said to be an {\tt attracting set} for the
semigroup $\{ T(h), h\ge 0\}$ on ${\cal K}^{+}$ in the topology
$\Theta^{loc}_{+}$ if for any bounded in ${\cal F}^b_{+}$ set
$B\subset {\cal K}^{+}$ one has $T(h)B\to P$ in the topology
$\Theta^{loc}_{+}$ when $h\to +\infty.$ }

\medskip

\noindent {\bf Definition~7.} {\it A set ${\cal U} \subset {\cal
K}^{+}$ is said to be a {\tt trajectory attractor} of semigroup
$\{ T(h), h\ge 0\}$ if
\begin{itemize}
\item[1)] ${\cal U}$ is bounded in ${\cal F}^b_{+}$ and compact in
$\Theta^{loc}_{+};$
 \item[2)] ${\cal U}$ is strictly invariant under the semigroup $\{ T(h),
h\ge 0\}$ i.e., $T(h){\cal U}={\cal U}$ for all $h\ge 0;$
\item[3)] ${\cal U}$ is an attracting set for the semigroup $\{
T(h), h\ge 0\}$ on ${\cal K}^{+}$ in the topology
$\Theta^{loc}_{+}.$

\end{itemize} }

\medskip

\noindent {\bf Theorem~5.} {\it Under the assumptions of Theorem~4
the semigroup $\{ T(h), h\ge 0\}$ on ${\cal K}^{+}$ possesses the
trajectory attractor ${\cal U}$.}

\smallskip

\noindent {\it Proof of Theorem~5.} Since any weak limiting
solution of (\ref{sdd3-1}), (\ref{sdd3-ic}) is a *-weak limit of
Galerkin approximate solutions to (\ref{sdd3-g}), (\ref{sdd3-ic})
we deduce some estimates for these approximate solutions.

Using (\ref{sdd3-F1}) and (\ref{sdd3-11}), we have for
$u(t)=u^{m,n}(t)$ :
$${d\over dt}\Vert u(t)\Vert^2 + 2 \Vert A^{1/2}u(t)\Vert^2 + 2d
\Vert u(t)\Vert^2 \le 2M_f|\Omega |^{3/2} C_b C_{\xi,1}\cdot \Vert
u(t)\Vert \le d \Vert u(t)\Vert^2 + \frac{M^2_f|\Omega |^{3} C^2_b
C^2_{\xi,1}}{d}.$$ This and $\Vert A^{1/2}u(t)\Vert^2 \ge
\lambda_1\Vert u(t)\Vert^2$ give
$$ {d\over dt}\Vert u(t)\Vert^2 + 2(d+\lambda_1) \Vert u(t)\Vert^2\le
\frac{M^2_f|\Omega |^{3} C^2_b C^2_{\xi,1}}{d}.
$$
Multiplying by $e^{(d+2\lambda_1)t}$ and integrating from $0$ to
$t,$ one obtains
\begin{equation}\label{sdd3-18}
\Vert u(t)\Vert^2 \le \Vert u(0)\Vert^2 \cdot e^{-2(d+\lambda_1)t}
+\frac{M^2_f|\Omega |^{3} C^2_b C^2_{\xi,1}}{d(d+2\lambda_1)},
\quad t\ge 0.
\end{equation}

Now we  multiply (\ref{sdd3-6}) by $\dot g_{k,m}(t)$ and take the
sum over $k=1,..,m$, and then we multiply (\ref{sdd3-6}) by $
g_{k,m}(t)$ and take the sum again over $k=1,..,m.$ The sum of the
obtained equations is (for $u=u^m$)
$$
\begin{array}{ll}
&\quad
 \langle F_n(u_t),\dot u(t)+u(t)\rangle \\
 &={1\over 2} {d\over dt} \left\{ \Vert A^{1\over 2}
u(t)\Vert^2+ (d+1)\Vert
 u(t)\Vert^2 \right\} + \Vert\dot u(t)\Vert^2 +
 \Vert A^{1\over 2} u(t)\Vert^2+d\Vert u(t)\Vert^2.
 \end{array}
 $$
Using (\ref{sdd3-F1}), we obtain positive constants $\gamma_1,
d_1$ (independent of $m$ and $n$) such that
\begin{equation}\label{sdd3-19}
 {d\over dt} \Psi (t) + \gamma_1\Psi (t) \le d_1, \quad
 \hbox{ where } \quad \Psi (t)\equiv \Vert A^{1\over 2}
u(t)\Vert^2+ (d+1)\Vert u(t)\Vert^2.
 \end{equation}
Multiplying it by $e^{\gamma_1 t}$ and integrating from $\tau>0$
to $\tau +h$ ($h>0$), we get $\Psi (\tau+h)e^{\gamma_1 (\tau+h)}
\le \Psi (\tau)e^{\gamma_1 (\tau)}+ d_1\gamma_1^{-1}e^{\gamma_1
(\tau+h)}.$ It gives $\Psi (\tau+h)\le \Psi (\tau)e^{-\gamma_1 h}+
d_1\gamma_1^{-1}.$ Integrating from $\tau=0$ to $\tau=1,$ one gets
$$\begin{array}{cl} \int^1_0\Psi (\tau+h) d\tau =
\int^{h+1}_h\Psi (s) ds \le & e^{-\gamma_1 h} \cdot \int^{1}_0\Psi
(s) ds
 + d_1\gamma_1^{-1}.
\end{array}
$$
Using the last inequality and definition of $\Psi$ (see
(\ref{sdd3-19})), we obtain
$$\int^{h+1}_h (\Vert A^{1\over 2} u(s)\Vert^2+ (d+1)\Vert
u(s)\Vert^2) ds \le e^{-\gamma_1 h} \cdot \int^{1}_0(\Vert
A^{1\over 2} u(s)\Vert^2+ (d+1)\Vert u(s)\Vert^2) ds
 + d_1\gamma_1^{-1}.
$$
 This and estimate (\ref{sdd3-14}) give that
\begin{equation}\label{sdd3-20}
\int^{h+1}_h \Vert A^{1\over 2} u(s)\Vert^2 ds \le e^{-\gamma_1 h}
\cdot \left( d+{3\over 2}\right) \left[\left( \Vert
u(0)\Vert^2+\tilde k_3\right)e^{\tilde k_1} - \tilde k_3\right] +
d_1\gamma_1^{-1}.
\end{equation}
It is important to note that all the constants $\gamma_1, d_1,
\tilde k_1, \tilde k_3$ are independent of $m$ and $n.$

In the same way, properties (\ref{sdd3-18}), (\ref{sdd3-20}) allow
one to get from (\ref{sdd3-6}) the following estimate
\begin{equation}\label{sdd3-21}
\int^{h+1}_h \Vert A^{-{1/2}} \dot u(s)\Vert^2 ds \le e^{-\gamma_2
h} \cdot d_2 \left( \Vert u(0)\Vert^2+1 \right)+ d_3
\end{equation}
with positive constants $\gamma_2, d_2, d_3$ independent of $m$
and $n.$

Estimates (\ref{sdd3-18}), (\ref{sdd3-20}) and (\ref{sdd3-21})
imply that any approximate solution $u^m=u^{n,m}$ belongs to
${\cal F}^b_{+}$. Moreover, there exists $R_1>0,$ such that the
ball $B_{R_1}=\{ w(\cdot)\in {\cal F}^b_{+} \quad | \quad
||w||_{{\cal F}^b_{+}} \le R_1 \}$ is an absorbing set for all the
approximate solutions $u^{m}=u^{n,m}$ of the problem
(\ref{sdd3-g}), (\ref{sdd3-ic}).
\medskip

 The constant $R_1>0$ is independent of $m$
and $n$ which gives that the ball $B_{R_1}=\{ w(\cdot)\in {\cal
F}^b_{+} \quad | \quad ||w||_{{\cal F}^b_{+}} \le R_1 \}$  is
also absorbing for any *-weak limit in the space 
$X_T$ (see (\ref{sdd3-13}), (\ref{sdd3-*})) of a subsequence $\{
u^{n_k,m_k} \} \subset \{ u^{n,m} \}.$ Particularly, this ball is
absorbing for any weak solution of (\ref{sdd3-g}), (\ref{sdd3-ic})
and for any  weak limiting solution of (\ref{sdd3-1}),
(\ref{sdd3-ic}). Hence it is an attracting (in the topology
$\Theta^{loc}_{+}$) set and by Remark~5 it is compact in
$\Theta^{loc}_{+}$. These properties together with Remarks 6,7
allow us to apply Theorem~3.1 from
\cite{Chepyzhov-Vishik_JMPA-1997} which completes the proof of
Theorem~5. \rule{5pt}{5pt}
\medskip

As an application we can consider the diffusive Nicholson's
blowflies equation (see e.g. So and Yang (1998), So, Wu and Yang
(2000)) with state-dependent (both discrete and distributed)
delays. More precisely, we consider equations (\ref{sdd3-1}) and
(\ref{sdd3-g}) where $-A$ is the Laplace operator with the
Dirichlet boundary conditions, $\Omega\subset R^{n_0}$ is a
bounded domain with a smooth boundary, the function $f$ can be a
constant as in So and Yang (1998), So, Wu and Yang (2000) which
leads to the local in space coordinate term or, for example, $
f(s)={1\over \sqrt{4\pi\alpha}} e^{-s^2/4\alpha}$, as in So, Wu,
Zou (2001) which corresponds to the non-local term, the nonlinear
function $b$ is given by $b(w)=p\cdot we^{-w}.$ Function $b$ is
bounded, so the conditions of Theorems~1-4 are satisfied. As a
result, we conclude that for any functions $\xi^n$ satisfying
conditions of
Theorems~1 and 2 
the dynamical system $(S_t,H)$ has a global attractor (Theorem~3).
Assuming the discrete delay $\eta$ is locally Lipschitz we get
(Theorem~4) the existence of weak limiting solutions of
(\ref{sdd3-1}), (\ref{sdd3-ic}) and the existence of trajectory
attractor (Theorem~5) for the corresponding translation semigroup.

\medskip

\noindent {\bf Acknowledgements.} The work was supported in part
by INTAS. The author wishes to thank Professor Hans-Otto Walther
for bringing state-dependent delay differential equations to his
attention during the visit to the Institute of Mathematics at
University of Giessen in 2002.

\bigskip
\hfill March 22, 2005

\end{document}